\documentclass[12pt]{amsart}

\usepackage{amssymb}
\usepackage{graphicx}
\usepackage{hyperref}
\usepackage{ifthen}
\usepackage{pict2e}
\usepackage{xargs}
\usepackage{xspace}

\newcommand{\definedterm}[1]{\emph{#1}}

\newcommand{\AD}{\mathtt{AD}}
\newcommand{\ADR}{\mathtt{AD}_\R}
\newcommand{\additivity}[1]{\mathtt{add}(#1)}
\newcommand{\auxiliaryfamily}[2]{\calF^{#1}_{\mspace{-0.6mu} #2}}
\newcommand{\auxiliarygraph}[1]{#1^*}
\newcommand{\Bairespace}[1][]{
  \ifthenelse{\equal{#1}{}}{\functions{\N}{\N}}{\functions{#1}{\N}}
}
\newcommand{\Borelchromaticnumber}[1]{\chi_{\text{B}}(#1)}
\newcommand{\bbE}{\mathbb{E}}
\newcommand{\bbG}{\mathbb{G}}
\newcommand{\bbN}{\mathbb{N}}
\newcommand{\bbR}{\mathbb{R}}
\newcommand{\bbZ}{\mathbb{Z}}
\newcommand{\calF}{\mathcal{F}}
\newcommand{\calM}{\mathcal{M}}
\newcommand{\calN}{\mathcal{N}}
\newcommand{\Cantorspace}[1][]{
  \ifthenelse{\equal{#1}{}}{\functions{\N}{2}}{\functions{#1}{2}}
}
\newcommand{\Cantortree}[1][]{
  \ifthenelse{\equal{#1}{}}{\functions{<\N}{2}}{\functions{#1}{2}}
}
\newcommand{\cardinality}[1]{|#1|}
\newcommand{\chromaticnumber}[1]{\chi(#1)}
\newcommand{\closedinterval}[3][]{[#2, #3]_{#1}}
\newcommand{\comparable}[1]{\parallel_{#1}}
\newcommand{\composition}{\circ}
\newcommandx{\concatenation}[2][1 = undefined, 2 = undefined]{
  \ifthenelse{\equal{#1}{undefined}}{{}\smallfrown}{
    \ifthenelse{\equal{#2}{undefined}}{\smallfrown_{#1}}{\bigoplus_{#1} #2}
  }
}

\newcommand{\DC}{\mathtt{DC}}
\newcommandx{\Deltaclass}[2][1=,2=]{
  \ifthenelse{\equal{#2}{}}{\mathbf{\Delta}_{#1}}{\mathbf{\Delta}^{#1}_{#2}}
}
\newcommand{\derivative}[2][]{
  \ifthenelse{\equal{#1}{}}{#2'}{#2^{(#1)}}
}
\newcommand{\emptysequence}{\emptyset}
\newcommand{\equivalenceclass}[2]{[#1]_{#2}}
\newcommand{\equivalencerelation}[1]{\equiv_{#1}}
\newcommand{\existinfinitelymany}{\exists^\infty}
\newcommand{\extendedby}{\sqsubseteq}
\newcommand{\extensions}[1]{\calN_{#1}}
\newcommand{\Ezero}{\bbE_0}
\newcommand{\finitechromaticnumber}[1]{\chi_{\text{fin}}(#1)}
\newcommand{\finitesubsets}[1]{\sets{< \aleph_0}{#1}}
\newcommand{\forcomeagerlymany}{\forall^*}
\newcommand{\from}{\colon}
\newcommandx{\functions}[3][3 =]{
  \ifthenelse{\equal{#3}{}}{#2^{#1}}{#2^{#1}_{#3}}
}

\newcommand{\Gzero}{\bbG_0}
\newcommand{\horizontalsection}[2]{#1^{#2}}
\newcommand{\image}[2]{#1(#2)}
\newcommand{\incomparable}[1]{\perp_{#1}}
\newcommandx{\intersection}[2][1 =, 2 =]{
  \ifthenelse{\equal{#1}{}}{\cap}{
    \ifthenelse{\equal{#2}{}}{\bigcap #1}{{\bigcap_{#1} #2}}
  }
}
\newcommand{\into}{\hookrightarrow}
\newcommand{\length}[1]{|#1|}
\newcommand{\mathcomma}{\text{, }}
\newcommand{\mathcommaand}{\text{, and }}
\newcommand{\mathor}{\text{ or }}
\newcommand{\meagerideal}{\calM}
\newcommand{\N}{\bbN}
\newcommand{\openclosedinterval}[3][]{(#2, #3]_{#1}}
\newcommand{\pair}[2]{(#1, #2)}
\newcommandx{\Piclass}[2][1=,2=]{
  \ifthenelse{\equal{#2}{}}{\mathbf{\Pi}_{#1}}{\mathbf{\Pi}^{#1}_{#2}}
}
\newcommand{\positiveintegers}{\Z^+}
\newcommand{\predecessorordinal}[1]{\boldsymbol{\lambda}^1_{#1}}
\newcommand{\preimage}[2]{#1^{-1}(#2)}
\newcommand{\projection}[1][]{
  \ifthenelse{\equal{#1}{}}{\mathrm{proj}}{\mathrm{proj}_{#1}}
}
\newcommand{\projectiveordinal}[1]{\boldsymbol{\delta}^1_{#1}}
\newcommandx{\rank}[2][2 = ]{
  \ifthenelse{\equal{#2}{}}{\rho(#1)}{\rho_{#1}(#2)}
}
\renewcommand{\restriction}[2]{#1 \upharpoonright #2}
\newcommandx{\sequence}[2][2 = undefined]{
  \ifthenelse{\equal{#2}{undefined}}{(#1)}{
    (#1)_{#2}
  }
}
\newcommandx{\set}[2][2 = undefined]{
  \ifthenelse{\equal{#2}{undefined}}{\{ #1 \}}{
    \{ #1 \suchthat #2 \}
  }
}
\newcommand{\setcomplement}[1]{\twiddle #1}
\newcommand{\setcomplementsubscript}[1]{\twiddlesubscript #1}
\newcommand{\sets}[2]{\functions{#1}{[#2]}}
\newcommandx{\Sigmaclass}[2][1=,2=]{
  \ifthenelse{\equal{#2}{}}
    {\mathbf{\Sigma}_{#1}}
    {\mathbf{\Sigma}^{#1}_{#2}}
}
\newcommand{\strictquasiorder}[1]{<_{#1}}
\newcommand{\successor}[1]{#1^+}
\newcommand{\suchthat}{\mid}
\newcommand{\symmetrization}[1]{#1^{\pm 1}}
\newcommand{\R}{\bbR}
\newcommand{\topology}[1]{\tau_{#1}}
\newcommand{\twiddle}
  {\raisebox{1.5pt}{\scalebox{.75}{$\mathord{\sim}$}}}
\newcommand{\twiddlesubscript}
  {\raisebox{1pt}{\scalebox{.5}{$\mathord{\sim}$}}}
\newcommandx{\union}[2][1 =, 2 =]{
  \ifthenelse{\equal{#1}{}}{\cup}{
    \ifthenelse{\equal{#2}{}}{\bigcup #1}{{\bigcup_{#1} #2}}
  }
}
\newcommand{\verticalsection}[2]{#1_{#2}}
\newcommand{\Z}{\bbZ}
\newcommand{\ZF}{\mathtt{ZF}}

\newcommand{\Baire}{Baire\xspace}
\newcommand{\Borel}{Bor\-el\xspace}
\newcommand{\Dilworth}{Dil\-worth\xspace}
\newcommand{\Hausdorff}{Haus\-dorff\xspace}

\newcommand{\Kada}{Ka\-da\xspace}
\newcommand{\Kanovei}{Kan\-o\-vei\xspace}
\newcommand{\Kuratowski}{Kur\-at\-ow\-ski\xspace}

\newcommand{\Montgomery}{Mont\-gom\-er\-y\xspace}
\newcommand{\Novikov}{No\-vik\-ov\xspace}

\newcommand{\Souslin}{Sous\-lin\xspace}
\newcommand{\Ulam}{U\-lam\xspace}
\newcommand{\Woodin}{Wood\-in\xspace}

\newenvironment{lemmaproof}{
  
  \begin{proof}
}{\end{proof}}

\newenvironment{propositionproof}{
  
  \begin{proof}
}{\end{proof}}

\newtheorem{lemma}{Lemma}
\newtheorem{proposition}[lemma]{Proposition}

\newtheorem{theorem}[lemma]{Theorem}

\theoremstyle{definition}
\newtheorem{remark}[lemma]{Remark}

\begin{document}


\begin{abstract}
  We generalize \Kada's definable strengthening of \Dilworth's
  characterization of the class of quasi-orders admitting an antichain of
  a given finite cardinality. 
\end{abstract}

\author[R. Carroy]{Rapha\"{e}l Carroy}

\address{
  Rapha\"{e}l Carroy \\
  Dipartimento di Matematica ``Giuseppe Peano'' \\
  Universit\`{a} di Torino \\
  Palazzo Campana \\
  Via Carlo Alberto 10 \\
  10123 Torino, Italia
}

\email{raphael.carroy@unito.it}

\urladdr{
  \url{http://www.logique.jussieu.fr/~carroy/indexeng.html}
}

\author[B.D. Miller]{Benjamin D. Miller}

\address{
  Benjamin D. Miller \\
  Kurt G\"{o}del Research Center for Mathematical Logic \\
  Universit\"{a}t Wien \\
  Augasse 2--6 \\
  1090 Wien \\
  Austria
 }

\email{benjamin.miller@univie.ac.at}

\urladdr{
  \url{http://www.logic.univie.ac.at/benjamin.miller}
}

\author[Z. Vidny\'{a}nszky]{Zolt\'{a}n Vidny\'{a}nszky}

\address{
  Zolt\'{a}n Vidny\'{a}nszky \\
  Kurt G\"{o}del Research Center for Mathematical Logic \\
  Universit\"{a}t Wien \\
  Augasse 2--6 \\
  1090 Wien \\
  Austria
 }

\email{zoltan.vidnyanszky@univie.ac.at}

\urladdr{
  \url{http://www.logic.univie.ac.at/~vidnyanszz77}
}



\title[The existence of small antichains]{On the existence of small
  antichains for definable quasi-orders}

\maketitle

\section*{Introduction}

A binary relation $R$ on a set $X$ is a \definedterm{quasi-order} if it
is reflexive and transitive. Two points $x, y \in X$ are \definedterm
{$R$-comparable} if $x \mathrel{R} y$ or $y \mathrel{R} x$, and
\definedterm{$R$-incomparable} otherwise. A set $Y \subseteq
X$ is an \definedterm{$R$-chain} if any two points of $Y$ are
$R$-comparable, and an \definedterm{$R$-antichain} if any two
distinct points of $Y$ are $R$-incomparable.

\Dilworth showed that if $k \in \positiveintegers$, $X$ is finite, and there
is no $R$-antichain of cardinality $k + 1$, then there is a cover
$\sequence{C_i}[i < k]$ of $X$ by $R$-chains (see \cite[Theorem 1.1]
{Dilworth}).

A subset of a topological space $X$ is \definedterm{\Borel} if it is in the
$\sigma$-algebra generated by the topology $\topology{X}$ of $X$,
\definedterm{analytic} if it is a continuous image of a closed subset of
$\Bairespace$, and \definedterm{$\aleph_0$-universally \Baire} if its
preimage under any continuous function $\phi \from \Cantorspace \to
X$ has the \Baire property.

Here we establish the following strengthening of \Dilworth's theorem:

\begin{theorem} \label{main}
  Suppose that $k \in \positiveintegers$, $X$ is a \Hausdorff space, and
  $R$ is an $\aleph_0$-universally-\Baire quasi-order on $X$ whose
  incomparability relation is analytic. Then exactly one of the following
  holds:
  \begin{enumerate}
    \item There is a cover $\sequence{C_i}[i < k]$ of $X$ by \Borel
      $R$-chains.
    \item There is an $R$-antichain of cardinality $k + 1$.
  \end{enumerate}
\end{theorem}

The \definedterm{equivalence relation} on $X$ associated with $R$ is
that with respect to which two points $x, y \in X$ are equivalent if $x
\mathrel{R} y$ and $y \mathrel{R} x$, and the \definedterm{strict
relation} associated with $R$ is that with respect to which two
points $x, y \in X$ are related if $x \mathrel{R} y$ but $\neg y \mathrel{R} x$.
\Kada established the special case of Theorem \ref{main} in which the
strict quasi-order is co-analytic and both the equivalence and incomparability
relations are analytic (see \cite[Theorem $1'$]{Kada}). Whereas his
intricate argument relied heavily upon recursion-theoretic methods, we
utilize only elementary ideas and the $\Gzero$ dichotomy (see \cite
[Theorem 6.3]{KST}), which itself has a classical proof (see \cite
[Theorem 8]{Miller}).

A subset of an analytic \Hausdorff space is \definedterm{$\Sigmaclass
[1][1]$} if it is analytic. More generally, for each $n \in \positiveintegers$,
a subset of an analytic \Hausdorff space is \definedterm{$\Piclass[1][n]$}
if its complement is $\Sigmaclass[1][n]$, and \definedterm{$\Sigmaclass
[1][n+1]$} if it is a continuous image of a $\Piclass[1][n]$ subset of an
analytic \Hausdorff space. A subset of an analytic \Hausdorff space
is \definedterm{$\Deltaclass[1][n]$} if it is both $\Sigmaclass[1][n]$ and
$\Piclass[1][n]$. \Souslin's theorem ensures that the families of \Borel
and $\Deltaclass[1][1]$ sets coincide (see, for example, \cite[Theorem
28.1]{Kechris}). The axiom of determinacy ($\AD$) implies that the
family of $\Deltaclass[1][2n+1]$ sets has a rich structural theory
analogous to that of the \Borel sets (see, for example, \cite{Jackson}).

We also obtain the following analog of Theorem \ref{main} under
determinacy:

\begin{theorem}[$\AD$] \label{main:AD}
  Suppose that $k \in \positiveintegers$, $n \in \N$, $X$ is an analytic
  \Hausdorff space, and $R$ is a quasi-order on $X$ whose
  incomparability relation is $\Sigmaclass[1][2n+1]$. Then exactly one
  of the following holds:
  \begin{enumerate}
    \item There is a cover $\sequence{C_i}[i < k]$ of $X$ by $\Deltaclass
      [1][2n+1]$ $R$-chains.
    \item There is an $R$-antichain of cardinality $k + 1$.
  \end{enumerate}
\end{theorem}

In addition, we generalize \Dilworth's theorem to arbitrary quasi-orders
on analytic \Hausdorff spaces under the strengthening of determinacy
where the players specify elements of $\R$ instead of $\N$ ($\ADR$):

\begin{theorem}[$\ADR$] \label{main:ADR} 
  Suppose that $k \in \positiveintegers$, $X$ is an analytic \Hausdorff
  space, and $R$ is a quasi-order on $X$. Then exactly one of the
  following holds:
  \begin{enumerate}
    \item There is a cover $\sequence{C_i}[i < k]$ of $X$ by $R$-chains.
    \item There is an $R$-antichain of cardinality $k + 1$.
  \end{enumerate}
\end{theorem}

In \S1, we establish Theorem \ref{main}. In \S2, we describe the minor
alterations to the proof necessary to obtain Theorems \ref{main:AD}
and \ref{main:ADR}. We work in the base theory $\ZF + \DC$
throughout.

\section{The classical case}

A binary relation $G$ on a set $X$ is a \definedterm{graph} if it is
irreflexive and symmetric. A ($Y$-)\definedterm{coloring} of $G$ is a
function $c \from X \to Y$ such that $w \mathrel{G} x \implies c(w)
\neq c(x)$ for all $w, x \in X$. The \definedterm{chromatic number} of
$G$, written $\chromaticnumber{G}$, is the least cardinal $\kappa$
for which there is a $\kappa$-coloring of $G$ (if such a cardinal
exists). We use $\finitechromaticnumber{G}$ to denote the supremum
of the chromatic numbers of the graphs of the form $\restriction{G}
{F}$, where $F \subseteq X$ is a finite set. We use $\auxiliarygraph
{G}$ to denote the supergraph of $G$ with respect to which two points
$x, y \in X$ are related if and only if there is a finite superset $F
\subseteq X$ of $\set{x, y}$ such that $c(x) \neq c(y)$ for every
$\finitechromaticnumber{G}$-coloring $c$ of $\restriction{G}{F}$.
Note that if $\finitechromaticnumber{G} = \aleph_0$, then $G =
\auxiliarygraph{G}$.

\begin{proposition} \label{union}
  Suppose that $X$ is a set, $G$ is a graph on $X$, and $G' \subseteq
  \auxiliarygraph{G}$ is finite. Then there is a finite set $F \subseteq X$
  containing $\union[i < 2][\image{\projection[i]}{G'}]$ such that every
  $\finitechromaticnumber{G}$-coloring $c$ of $\restriction{G}{F}$ is a
  coloring of $\symmetrization{(G')}$.
\end{proposition}

\begin{propositionproof}
  For all $\pair{x}{y} \in G'$, fix a finite superset $F_{\pair{x}{y}}
  \subseteq X$ of $\set{x, y}$ such that $c(x) \neq c(y)$ for every
  $\finitechromaticnumber{G}$-coloring $c$ of $\restriction{G}{F_{\pair
  {x}{y}}}$, and observe that the set $F = \union[\pair{x}{y} \in G']
  [F_{\pair{x}{y}}]$ is as desired.
\end{propositionproof}

A set $Y \subseteq X$ is a \definedterm
{$G$-clique} if any two distinct points of $Y$ are $G$-related, and
\definedterm{$G$-independent} if no two points of $Y$ are $G$-related.

\begin{proposition} \label{clique}
  Suppose that $X$ is a set, $G$ is a graph on $X$, and $C \subseteq
  X$ is a finite $\auxiliarygraph{G}$-clique. Then $\cardinality{C} \le
  \finitechromaticnumber{G}$.
\end{proposition}

\begin{propositionproof}
  By Proposition \ref{union}, there is a finite set $F \subseteq X$
  containing $C$ such that $\restriction{c}{C}$ is injective for every
  $\finitechromaticnumber{G}$-coloring $c$ of $\restriction{G}{F}$, in
  which case the pigeon-hole principle ensures that $\cardinality{C} \le
  \finitechromaticnumber{G}$.
\end{propositionproof}

The \definedterm{horizontal sections} of a set $R \subseteq X \times Y$
are the sets of the form $\horizontalsection{R}{y} = \set{x \in X}[x
\mathrel{R} y]$, where $y \in Y$. The \definedterm{vertical sections} are
the sets of the form $\verticalsection{R}{x} = \set{y \in Y}[x \mathrel{R}
y]$, where $x \in X$.

\begin{proposition} \label{antichain}
  Suppose that $X$ is a set, $G$ is a graph on $X$ for which
  $\finitechromaticnumber{G} < \aleph_0$, $x, y \in X$, and there is a
  $\auxiliarygraph{G}$-clique $C \subseteq \verticalsection
  {\auxiliarygraph{G}}{x} \union \verticalsection{\auxiliarygraph{G}}{y}$
  of cardinality $\finitechromaticnumber{G}$. Then $x \mathrel
  {\auxiliarygraph{G}} y$.
\end{proposition}

\begin{propositionproof}
  Proposition \ref{union} yields a finite set $F \subseteq X$ containing
  $C \union \set{x, y}$ such that $\restriction{c}{C}$ is injective and
  $\forall w \in \set{x, y} \forall z \in C \intersection \verticalsection
  {\auxiliarygraph{G}}{w} \ c(w) \neq c(z)$ for every
  $\finitechromaticnumber{G}$-coloring $c$ of $\restriction{G}{F}$. But
  if $c$ is such a coloring, then $\image{c}{C} = \finitechromaticnumber
  {G}$, so $c(x) \in \image{c}{C \intersection \verticalsection
  {\auxiliarygraph{G}}{y}}$, thus $c(x) \neq c(y)$, hence $x \mathrel
  {\auxiliarygraph{G}} y$.
\end{propositionproof}

We use $\comparable{R}$, $\equivalencerelation{R}$, $\incomparable
{R}$, and $\strictquasiorder{R}$ to denote the comparability,
equivalence, incomparability, and strict relations associated with $R$.

\begin{proposition} \label{transitive}
  Suppose that $X$ is a set and $R$ is a quasi-order on $X$. Then $R
  \setminus \mathord{\auxiliarygraph{\incomparable{R}}}$ is transitive.
\end{proposition}

\begin{propositionproof}
  Suppose, towards a contradiction, that there exist $x, y, z \in X$ for
  which $x \mathrel{(R \setminus \mathord{\auxiliarygraph
  {\incomparable{R}}})} y \mathrel{(R \setminus \mathord{\auxiliarygraph
  {\incomparable{R}}})} z$, as well as a finite set $F \subseteq X$
  containing $\set{x, z}$ such that $c(x) \neq c(z)$ for every
  $\finitechromaticnumber{\incomparable{R}}$-coloring $c$ of
  $\restriction{\mathord{\incomparable{R}}}{F}$. Then $x \mathrel{R} z$,
  so $x$ and $z$ are not $\incomparable{R}$-related, thus
  $\finitechromaticnumber{\incomparable{R}} < \aleph_0$. For all $w \in
  \set{x, z}$, the fact that $w$ and $y$ are not $\auxiliarygraph
  {\incomparable{R}}$-related yields an $R$-chain $C_w \subseteq F
  \union \set{y}$ containing $\set{w,y}$ for which $(F \union \set{y})
  \setminus C_w$ is a union of $\finitechromaticnumber{\incomparable
  {R}} - 1$ $R$-chains, and therefore does not contain an $R$-antichain
  of cardinality $\finitechromaticnumber{\incomparable{R}}$. Then the
  set $C = (C_x \intersection \horizontalsection{R}{y}) \union (C_z
  \intersection \verticalsection{R}{y})$ is an $R$-chain containing $\set
  {x,z}$, so $(F \union \set{y}) \setminus C$ is not a union of
  $\finitechromaticnumber{\incomparable{R}} - 1$ $R$-chains, thus
  \Dilworth's theorem yields an $R$-antichain $A \subseteq (F \union
  \set{y}) \setminus C$ of cardinality $\finitechromaticnumber
  {\incomparable{R}}$. Fix $u \in A \intersection C_x$ and $w \in A
  \intersection C_z$. As $u, w \notin C$, it follows that neither $u \mathrel
  {R} y$ nor $y \mathrel{R} w$, so the fact that $C_x$ and $C_z$ are
  $R$-chains ensures that $w \strictquasiorder{R} y \strictquasiorder{R}
  u$, contradicting the fact that $A$ is an $R$-antichain.
\end{propositionproof}

Define $\closedinterval[R]{x}{y} = \set{z \in X}[x \mathrel{R} z \mathrel
{R} y]$ and $\openclosedinterval[R]{x}{y} = \closedinterval[R]{x}{y}
\setminus \equivalenceclass{x}{\equivalencerelation{R}}$. We use
$\mathord{\concatenation}$, $\extendedby$, and $\sequence{i}$ to
denote concatenation, extension, and the sequence of length one
whose sole entry is $i$. Fix sequences $s_n \in \Cantorspace[n]$ that
are \definedterm{dense} in $\Cantortree$, in the sense that $\forall s \in
\Cantortree \exists n \in \N \ s \extendedby s_n$, and define $\Gzero =
\set{\pair{s_n \concatenation \sequence{i} \concatenation c}{s_n
\concatenation \sequence{1 - i} \concatenation c}}[c \in \Cantorspace
\mathcomma i < 2 \mathcommaand n \in \N]$.

\begin{proposition} \label{homomorphism}
  Suppose that $X$ is a topological space, $R$ is an
  $\aleph_0$-universally-\Baire quasi-order on $X$ that does not have
  antichains of arbitrarily large finite cardinality, and $\auxiliarygraph
  {\incomparable{R}}$ is $\aleph_0$-universally \Baire. Then there is no
  continuous homomorphism $\phi \from \Cantorspace \to X$ from
  $\Gzero$ to $\auxiliarygraph{\incomparable{R}}$.
\end{proposition}

\begin{propositionproof}
  As \Dilworth's theorem ensures that $\finitechromaticnumber
  {\incomparable{R}} < \aleph_0$, it is sufficient to show that if $\phi
  \from \Cantorspace \to X$ is a continuous homomorphism from
  $\Gzero$ to $\auxiliarygraph{\incomparable{R}}$, then there exists $x
  \in \image{\phi}{\Cantorspace}$ for which there is a continuous
  homomorphism from $\Gzero$ to $\restriction{\mathord{\auxiliarygraph
  {\incomparable{R}}}}{(\image{\phi}{\Cantorspace} \intersection
  \verticalsection{(\auxiliarygraph{\incomparable{R}})}{x})}$, since
  $\finitechromaticnumber{\incomparable{R}}$ applications of this fact
  yield a $\auxiliarygraph{\incomparable{R}}$-clique of cardinality
  $\finitechromaticnumber{\incomparable{R}} + 1$, contradicting
  Proposition \ref{clique}.
  
  Letting $G'$ be the pullback of $\auxiliarygraph{\incomparable{R}}$
  through $\phi \times \phi$, it is sufficient to find $c \in \Cantorspace$
  for which $\verticalsection{G'}{c}$ has the \Baire property and is not
  meager, as the proof of \cite[Proposition 6.2]{KST} ensures that every
  $\Gzero$-independent set with the \Baire property is meager, so \cite
  [Theorem 6.3]{KST} would then yield a continuous homomorphism
  $\psi \from \Cantorspace \to \verticalsection{G'}{c}$ from $\Gzero$ to
  $\restriction{\Gzero}{\verticalsection{G'}{c}}$ (although the existence
  of such a function also follows from a straightforward recursive
  construction), in which case the point $x = \phi(c)$ and the
  homomorphism $\phi \composition \psi$ are as desired.
  
  Suppose, towards a contradiction, that every vertical section of $G'$
  with the \Baire property is meager, and let $R'$ be the pullback of
  $R$ through $\phi \times \phi$. As $\auxiliarygraph{\incomparable
  {R}}$ and $R$ are $\aleph_0$-universally \Baire, the horizontal and
  vertical sections of $G'$ and $R'$ all have the \Baire property. As
  $\mathord{\incomparable{R'}} \subseteq G'$, every vertical section
  of $\incomparable{R'}$ is meager, and the \Kuratowski-\Ulam
  theorem (see, for example, \cite[Theorem 8.41]{Kechris}) ensures
  that $\comparable{R'}$ is comeager, so $R'$ is not meager.
    
  \begin{lemma} \label{homomorphism:nonmeager}
    There exists $\pair{b}{d} \in \Gzero$ for which $\closedinterval[R']{b}
    {d}$ is not meager.
  \end{lemma}
  
  \begin{lemmaproof}
    It is trivial to check that the binary relation $S'$ on $\Cantorspace$
    given by $c \mathrel{S'} d \iff \forcomeagerlymany b \in \Cantorspace
    \ (b \mathrel{R'} c \implies b \mathrel{R'} d)$ is a quasi-order, and for
    no $\pair{d}{c} \in \setcomplement{S'}$ is $\openclosedinterval[R']{c}
    {d}$ meager (see, for example, \cite[Proposition 8.26]{Kechris}). We
    can therefore assume that $\Gzero \subseteq S'$, so $\Gzero
    \subseteq \mathord{\equivalencerelation{S'}}$. As the smallest
    equivalence relation on $\Cantorspace$ containing $\Gzero$ is
    $\Ezero$ (by a straightforward inductive argument), it follows that
    $\Ezero \subseteq \mathord{\equivalencerelation{S'}}$. As the
    \Kuratowski-\Ulam and \Montgomery-\Novikov theorems (see, for
    example, \cite[Theorem 16.1]{Kechris}) ensure that for all $s \in
    \Cantortree$, the corresponding set $B_s = \set{c \in \Cantorspace}
    [\forcomeagerlymany b \in \extensions{s} \ b \mathrel{R'} c]$ has the
    \Baire property, and $c \equivalencerelation{S'} d \iff \forall s \in
    \Cantortree \ (c \in B_s \iff d \in B_s)$ for all $c, d \in \Cantorspace$,
    the fact that every $\Ezero$-invariant set with the \Baire property is
    meager or comeager (see, for example, \cite[Theorem 8.47]
    {Kechris}) yields a comeager $\equivalencerelation{S'}$-class.
    Fixing $s, t \in \Cantortree$ with the property that $R' \intersection
    (\extensions{s} \times \extensions{t})$ is comeager in $\extensions
    {s} \times \extensions{t}$, the \Kuratowski-\Ulam theorem implies
    that $\forcomeagerlymany c \in \extensions{t} \forcomeagerlymany
    b \in \extensions{s} \ b \mathrel{R'} c$, so $\forcomeagerlymany b,
    c \in \extensions{s} \ b \mathrel{R'} c$, thus there is an
    $\equivalencerelation{R'}$-class $C \subseteq \Cantorspace$ that
    is comeager in $\extensions{s}$. But non-meager subsets of
    $\Cantorspace$ with the \Baire property are not
    $\Gzero$-independent, and any pair $\pair{b}{d} \in \restriction
    {\Gzero}{C}$ is as desired.
  \end{lemmaproof}
  
  As $b \mathrel{G'} d$, Proposition \ref{transitive} ensures that $\forall
  c \in \closedinterval[R']{b}{d} \ (b \mathrel{G'} c \mathor c \mathrel{G'}
  d)$, so $\verticalsection{G'}{b}$ or $\verticalsection{G'}{d}$ is not
  meager, the desired contradiction.
\end{propositionproof}

\begin{remark}
  A similar approach can be used to eliminate the need for multiple
  applications of the $\Gzero$ dichotomy, and therefore the need to
  assume that $\additivity{\meagerideal} < \kappa$, in \cite{MV} (see
  \cite[Propositions 1.6.17 and 1.6.19]{Miller:Notes}).
\end{remark}

\begin{proposition} \label{independence}
  Suppose that $X$ is a set, $R$ is a quasi-order on $X$ that does not
  have antichains of arbitrarily large finite cardinality, $A \subseteq X$ is
  an $R$-antichain of cardinality $\finitechromaticnumber{\incomparable
  {R}}$, and $Y \subseteq X$ is $\auxiliarygraph{\incomparable
  {R}}$-independent. Then there exists $x \in A$ for which $\set{x}
  \union Y$ is $\auxiliarygraph{\incomparable{R}}$-independent.
\end{proposition}

\begin{propositionproof}
  Suppose, towards a contradiction, that there exists a function $\phi
  \from A \to Y$ whose graph is contained in $\auxiliarygraph
  {\incomparable{R}}$. As \Dilworth's theorem ensures that
  $\finitechromaticnumber{\incomparable{R}} < \aleph_0$, it follows
  that $A$ is a maximal $R$-antichain, and is therefore the union of the
  sets $A' = \set{x \in A}[A \intersection \horizontalsection{R}{\phi(x)}
  \neq \emptyset]$ and $A'' = \set{x \in A}[A \intersection
  \verticalsection{R}{\phi(x)} \neq \emptyset]$.
  
  \begin{lemma} \label{independence:disjoint}
    The sets $A'$ and $A''$ are disjoint.
  \end{lemma}
  
  \begin{lemmaproof}
    Suppose, towards a contradiction, that there exists $x \in A'
    \intersection A''$, and fix $y, z \in A$ for which $y \mathrel{R} \phi(x)
    \mathrel{R} z$. As $A$ is an $R$-antichain, it follows that $y = z$,
    so $\phi(x) \equivalencerelation{R} y$, thus the $\equivalencerelation
    {R}$-invariance of $\auxiliarygraph{\incomparable{R}}$ yields that
    $\phi(x) \mathrel{\auxiliarygraph{\incomparable{R}}} \phi(y)$,
    contradicting the $\auxiliarygraph{\incomparable{R}}$-independence
    of $Y$.
  \end{lemmaproof}
  
  \begin{lemma} \label{independence:A'}
    If $w', x' \in A'$ and $\phi(x') \mathrel{R} \phi(w')$, then $w' \mathrel
    {\auxiliarygraph{\incomparable{R}}} \phi(x')$.
  \end{lemma}
  
  \begin{lemmaproof}
    If $w'$ and $\phi(x')$ are not $\auxiliarygraph{\incomparable
    {R}}$-related, then $w' \comparable{R} \phi(x')$, so Lemma \ref
    {independence:disjoint} ensures that $w' \mathrel{(R \setminus
    \mathord{\auxiliarygraph{\incomparable{R}}})} \phi(x')$. But the
    $\auxiliarygraph{\incomparable{R}}$-independence of $Y$ implies
    that $\phi(x') \mathrel{(R \setminus \mathord{\auxiliarygraph
    {\incomparable{R}}})} \phi(w')$, thus Proposition \ref{transitive} yields
    that $w'$ and $\phi(w')$ are not $\auxiliarygraph{\incomparable
    {R}}$-related, a contradiction.
  \end{lemmaproof}
  
  \begin{lemma} \label{independence:A''}
    If $w'', x'' \in A''$ and $\phi(w'') \mathrel{R} \phi(x'')$, then $w'' \mathrel
    {\auxiliarygraph{\incomparable{R}}} \phi(x'')$.
  \end{lemma}
  
  \begin{lemmaproof}
    If $w''$ and $\phi(x'')$ are not $\auxiliarygraph{\incomparable
    {R}}$-related, then $w'' \comparable{R} \phi(x'')$, so Lemma \ref
    {independence:disjoint} ensures that $\phi(x'') \mathrel{(R \setminus
    \mathord{\auxiliarygraph{\incomparable{R}}})} w''$. But the
    $\auxiliarygraph{\incomparable{R}}$-independence of $Y$ implies
    that $\phi(w'') \mathrel{(R \setminus \mathord{\auxiliarygraph
    {\incomparable{R}}})} \phi(x'')$, thus Proposition \ref{transitive} yields
    that $\phi(w'')$ and $w''$ are not $\auxiliarygraph{\incomparable
    {R}}$-related, a contradiction.
  \end{lemmaproof}
  
  If $A' \neq \emptyset$, then the fact that $Y$ is an $R$-chain yields
  $x' \in A'$ for which $\phi(x')$ is $(\restriction{R}{\image{\phi}
  {A'}})$-minimal, so Lemma \ref{independence:A'} ensures that $A'
  \union \set{\phi(x')}$ is an $\auxiliarygraph{\incomparable{R}}$-clique,
  and since Lemma \ref{independence:disjoint} implies that $\phi(x')
  \notin A'$, Proposition \ref{clique} yields that $\cardinality{A'} <
  \finitechromaticnumber{\incomparable{R}}$. Similarly, if $A'' \neq
  \emptyset$, then the fact that $Y$ is an $R$-chain yields $x'' \in A''$
  for which $\phi(x'')$ is $(\restriction{R}{\image{\phi}{A''}})$-maximal,
  so Lemma \ref{independence:A''} ensures that $A'' \union \set{\phi
  (x'')}$ is an $\auxiliarygraph{\incomparable{R}}$-clique, and since
  Lemma \ref{independence:disjoint} implies that $\phi(x'') \notin A''$,
  Proposition \ref{clique} implies that $\cardinality{A''} <
  \finitechromaticnumber{\incomparable{R}}$. It follows that $A'$ and
  $A''$ are non-empty, so there are indeed $x' \in A'$ and $x'' \in A''$ for
  which $\phi(x')$ is $(\restriction{R}{\image{\phi}{A'}})$-minimal and
  $\phi(x'')$ is $(\restriction{R}{\image{\phi}{A''}})$-maximal. As $A
  \subseteq \verticalsection{(\auxiliarygraph{\incomparable{R}})}{\phi
  (x')} \union \verticalsection{(\auxiliarygraph{\incomparable{R}})}{\phi
  (x'')}$ by Lemmas \ref{independence:A'} and \ref{independence:A''},
  Proposition \ref{antichain} implies that $\phi(x') \mathrel
  {\auxiliarygraph{\incomparable{R}}} \phi(x'')$, contradicting the
  $\auxiliarygraph{\incomparable{R}}$-independence of $Y$.
\end{propositionproof}

For each $k \in \N$, let $\sets{k}{X}$ denote the family of all subsets
of $X$ of cardinality $k$, equipped with the topology generated by the
sets of the form $\set{F \in \sets{k}{X}}[\exists \pi \from F \into \calF
\ \forall x \in F \ x \in \pi(x)]$, where $\calF \in \sets{k}{\topology{X}}$.
Let $\sets{\le k}{X}$ denote the disjoint union of the spaces of the
form $\sets{j}{X}$, for $j \le k$. Similarly, let $\finitesubsets{X}$ denote
the disjoint union of the spaces of the form $\sets{k}{X}$, for $k \in \N$.
A set $Y \subseteq X$ \definedterm{punctures} a family $\calF
\subseteq \finitesubsets{X}$ if $F \intersection Y \neq \emptyset$ for
all $F \in \calF$.

\begin{proposition} \label{maximal}
  Suppose that $X$ is a \Hausdorff space, $G$ is an analytic graph on
  $X$ that admits a \Borel coloring $c \from X \to \N$, and $\calF
  \subseteq \finitesubsets{X}$ is an analytic set with the property that
  for every $G$-independent set $Y \subseteq X$, the corresponding
  set $\set{x \in X}[\set{x} \union Y \text{ is $G$-independent}]$
  punctures $\calF$. Then every $G$-independent \Borel subset of $X$
  is contained in a $G$-independent \Borel subset of $X$ that
  punctures $\calF$.
\end{proposition}

\begin{propositionproof}
  For each natural number $k$ and $G$-independent set $Y \subseteq
  X$, we use $\auxiliaryfamily{k}{Y}$ to denote the family of sets $F \in
  \calF$ with the property that $\cardinality{\set{x \in F}[\set{x} \union Y
  \text{ is not $G$-independent}]} \ge \cardinality{F} - k$. Note that
  $\auxiliaryfamily{0}{Y} = \emptyset$ and $\calF \intersection \sets{\le k}
  {X} \subseteq \auxiliaryfamily{k}{Y}$. It is sufficient to show that for all
  $k \in \N$, every $G$-independent \Borel set $B \subseteq X$ that
  punctures $\auxiliaryfamily{k}{B}$ is contained in a $G$-independent
  \Borel set $C \subseteq X$ that punctures $\auxiliaryfamily{k+1}{C}$,
  as repeated application of this fact yields an increasing sequence of
  $G$-independent \Borel supersets $B_k \subseteq X$ of any given
  $G$-independent \Borel subset of $X$ that puncture $\auxiliaryfamily
  {k}{B_k}$, in which case the set $\union[k \in \N][B_k]$ is as desired.
  
  Suppose that $k \in \N$, we have already established the
  aforementioned fact strictly below $k$, and $B \subseteq X$ is a
  $G$-independent \Borel set that punctures $\auxiliaryfamily{k}{B}$.
  Fix natural numbers $i_j$ such that $\forall i \in \N \existinfinitelymany
  j \in \N \ i = i_j$, and define $B_0' = B$. Given $j \in \N$ and a
  $G$-independent \Borel set $B_j' \subseteq X$ that punctures
  $\auxiliaryfamily{k}{B_j'}$, let $A_j'$ be the set of $x \in X$ for which
  there exists $F \in \calF$ disjoint from $B_j'$ with the property that $x
  \in F$ and $\cardinality{\set{y \in F \setminus \set{x}}[B_j' \union \set
  {y} \text{ is not $G$-independent}]} \ge \cardinality{F} - (k+1)$. The fact that
  $B_j'$ punctures $\auxiliaryfamily{k}{B_j'}$ ensures that $B_j' \union
  \set{x}$ is $G$-independent for all $x \in A_j'$, thus so too is $(A_j'
  \intersection \preimage{c}{\set{i_j}}) \union B_j'$. As the latter set is
  analytic, it is contained in a $G$-independent \Borel set (see, for
  example, the proof of \cite[Proposition 2]{Miller}), in which case $k$
  applications of the induction hypothesis yield a $G$-independent
  \Borel set $B_{j+1}' \subseteq X$ containing $(A_j' \intersection
  \preimage{c}{\set{i_j}}) \union B_j'$ that punctures $\auxiliaryfamily
  {k}{B_{j+1}'}$.
  
  To see that the $G$-independent \Borel set $C = \union[j \in \N][B_j']$
  punctures $\auxiliaryfamily{k+1}{C}$, observe that if $F \in
  \auxiliaryfamily{k+1}{C}$, then there exists $x \in F$ for which $C
  \union \set{x}$ is $G$-independent, as well as $j \in \N$ for which $F
  \in \auxiliaryfamily{k+1}{B_j'}$, and $j' \ge j$ for which $i_{j'} =
  c(x)$, in which case $B_{j'}' \intersection F \neq \emptyset$ or $x \in
  B_{j'+1}'$.
\end{propositionproof}

The \definedterm{\Borel chromatic number} of a graph $G$ on $X$ is
the least cardinal $\Borelchromaticnumber{G}$ of the form
$\cardinality{Y}$, where $Y$ is an analytic \Hausdorff space for which
there exists a \Borel $Y$-coloring of $G$ (if such a space exists).

\begin{proposition} \label{coloring}
  Suppose that $X$ is a \Hausdorff space and $R$ is a quasi-order on
  $X$ with the property that $\incomparable{R}$ is analytic and
  $\Borelchromaticnumber{\auxiliarygraph{\incomparable{R}}} \le
  \aleph_0$. Then $\Borelchromaticnumber{\auxiliarygraph
  {\incomparable{R}}} = \finitechromaticnumber{\incomparable{R}}$.
\end{proposition}

\begin{propositionproof}
  As the case $\finitechromaticnumber{\incomparable{R}} \in \set{1,
  \aleph_0}$ is trivial, suppose that $k \in \positiveintegers$, we have
  already established the proposition for $\finitechromaticnumber
  {\incomparable{R}} \le k$, and $\finitechromaticnumber{\incomparable
  {R}} = k + 1$. As $\auxiliarygraph{\incomparable{R}}$ is analytic,
  Propositions \ref{independence} and \ref{maximal} yield an
  $\auxiliarygraph{\incomparable{R}}$-independent \Borel set $B \subseteq
  X$ that intersects every $R$-antichain of cardinality $k + 1$. As
  \Dilworth's theorem ensures that $\finitechromaticnumber{\restriction
  {\mathord{\incomparable{R}}}{\setcomplement{B}}} = k$, the induction
  hypothesis yields a \Borel $k$-coloring $c$ of $\auxiliarygraph
  {(\restriction{\mathord{\incomparable{R}}}{\setcomplementsubscript
  {B}})}$. Observe that $\restriction{\mathord{\auxiliarygraph
  {\incomparable{R}}}}{\setcomplement{B}} \subseteq \auxiliarygraph
  {(\restriction{\mathord{\incomparable{R}}}{\setcomplementsubscript
  {B}})}$, for if $x, y \in \setcomplement{B}$ and $F \subseteq X$ is a
  finite set containing $\set{x, y}$ such that $d(x) \neq d(y)$ for every $(k+
  1)$-coloring $d$ of $\restriction{\mathord{\incomparable{R}}}{F}$,
  then $F \setminus B$ is a finite set containing $\set{x, y}$ such that
  $d(x) \neq d(y)$ for every $k$-coloring $d$ of $\restriction{\mathord
  {\incomparable{R}}}{(F \setminus B)}$. In particular, it follows that the
  extension of $c$ to $X$ with constant value $k$ on $B$ is a \Borel
  $(k+1)$-coloring of $\auxiliarygraph{\incomparable{R}}$.
\end{propositionproof}

As every analytic subset of a topological space is
$\aleph_0$-universally \Baire (see, for example, \cite[Theorem 21.6]
{Kechris}), Theorem \ref{main} follows from Proposition \ref
{homomorphism}, the $\Gzero$ dichotomy, and Proposition \ref
{coloring}.

\section{Generalizations under determinacy}

Given an ordinal $\alpha$, a subset of a topological space $X$ is
\definedterm{$\alpha$-\Borel} if it is in the closure of $\topology{X}$
under complements and unions of length strictly less than $\alpha$.
Given an aleph $\kappa$, a topological space is \definedterm
{$\kappa$-\Souslin} if it is a continuous image of a closed subset of
$\functions{\N}{\kappa}$.

For all $n > 0$, let $\projectiveordinal{n}$ denote the supremum of the
lengths of well-orders of the form $R / \mathord{\equivalencerelation
{R}}$, where $R$ is a $\Deltaclass[1][n]$ quasi-order on an analytic
\Hausdorff space. The axiom of determinacy ensures that the
$\Deltaclass[1][2n+1]$ and $\projectiveordinal{2n+1}$-\Borel subsets of
analytic \Hausdorff spaces coincide. It also yields an aleph
$\predecessorordinal{2n+1}$ for which $\projectiveordinal{2n+1} =
\successor{(\predecessorordinal{2n+1})}$, and implies that the
$\Sigmaclass[1][2n+1]$ and $\predecessorordinal{2n+1}$-\Souslin
subsets of analytic \Hausdorff spaces coincide (see, for example,
\cite{Jackson}).

A \definedterm{tree} on a set $I$ is a set $T \subseteq \functions{<\N}
{I}$ that is \definedterm{closed under initial segments}, in the sense
that $\forall t \in T \forall n < \length{t} \ \restriction{t}{n} \in T$. A
\definedterm{subtree} of $T$ is a tree $S \subseteq T$ on $I$. A
\definedterm{branch} through $T$ is a sequence $x \in \functions{\N}
{I}$ such that $\forall n \in \N \ \restriction{x}{n} \in T$. A tree is
\definedterm{well-founded} if it has no branches.

The \definedterm{pruning derivative} associates with each tree $T$ on
a set $I$ the subtree $\derivative{T} = \set{t \in T}[\exists i \in I \ t
\concatenation \sequence{i} \in T]$. The \definedterm{iterates} of the
pruning derivative are given by $\derivative[0]{T} = T$, $\derivative
[\alpha + 1]{T} = \derivative{(\derivative[\alpha]{T})}$ for all ordinals
$\alpha$, and $\derivative[\lambda]{T} = \intersection[\alpha <
\lambda][{\derivative[\alpha]{T}}]$ for all limit ordinals $\lambda$. The
\definedterm{pruning rank} of $T$ is the least ordinal $\rank{T}$ for
which $\derivative[\rank{T}]{T} = \derivative[\rank{T}+1]{T}$. A
straightforward induction shows that $T$ is well-founded if and only if
$\derivative[\rank{T}]{T} = \emptyset$. For each $t \in T$, let $\rank{T}
[t]$ denote the largest ordinal for which $t \in \derivative[{\rank{T}[t]}]
{T}$ (if such an ordinal exists).

An \definedterm{$(\alpha + 1)$-\Borel code} for a subset of $X$ is a pair
$\pair{f}{T}$, where $T$ is a well-founded tree on $\alpha \times \alpha$
and $f$ is a function associating to each sequence $t \in
\setcomplement{T}$ a subset of $X$ that is closed or open. Given such
a code, we recursively define $f^{(\beta)}$ on $\setcomplement
{\derivative[\beta]{T}}$ by setting $f^{(0)} = f$, letting $f^{(\beta+1)}$ be
the extension of $f^{(\beta)}$ given by $f^{(\beta+1)}(t) = \union[\gamma
< \alpha][{\intersection[\delta < \alpha][f^{(\beta)}(t \concatenation
\sequence{\pair{\gamma}{\delta}})]}]$ whenever $\rank{T}[t] = \beta$ for
all ordinals $\beta$, and defining $f^{(\lambda)} = \union[\beta <
\lambda][f^{(\beta)}]$ for all limit ordinals $\lambda$. The $(\alpha +
1)$-\Borel set \definedterm{coded} by $\pair{f}{T}$ is $f^{(\rank{T})}
(\emptysequence)$.

The proof of \Souslin's theorem shows that there is a function sending
each pair of functions witnessing that a set and its complement are
$\kappa$-\Souslin to a $(\kappa + 1)$-\Borel code for the set. Under
$\AD$, the coding lemma (see \cite[Lemma 7D.5]{Moschovakis}) and
projective uniformization (see, for example, \cite[Theorem 39.9]
{Kechris}) can be used to obtain a function sending each
$(\predecessorordinal{2n+1} + 1)$-\Borel code for a subset of an analytic
\Hausdorff space to a function witnessing that the encoded set is
$\predecessorordinal{2n+1}$-\Souslin.

\begin{proposition}[$\AD$] \label{maximal:AD}
  Suppose that $n \in \N$, $X$ is an analytic \Hausdorff space, $G$ is
  a $\Sigmaclass[1][2n+1]$ graph on $X$ that admits a $\Deltaclass[1]
  [2n+1]$ coloring $c \from X \to \predecessorordinal{2n+1}$, and
  $\calF \subseteq \finitesubsets{X}$ is a $\Sigmaclass[1][2n+1]$ set
  with the property that for every $G$-independent set $Y \subseteq
  X$, the corresponding set $\set{x \in X}[\set{x} \union Y \text{ is
  $G$-independent}]$ punctures $\calF$. Then every
  $G$-independent $\Deltaclass[1][2n+1]$ subset of $X$ is contained
  in a $G$-independent $\Deltaclass[1][2n+1]$ subset of $X$ that
  punctures $\calF$.
\end{proposition}

\begin{propositionproof}
  We proceed essentially as in the proof of Proposition \ref{maximal}.
  The first paragraph remains unchanged. The induction beginning
  in the second paragraph, however, has length $\predecessorordinal
  {2n+1}$ instead of $\omega$, which is problematic because naively
  applying \cite[Proposition 2]{Miller} at each stage of the induction
  requires too large a fragment of the axiom of choice. This problem
  can be alleviated by using the above remarks to keep track of codes
  for the sets $B_j'$ that are built along the way, which can be
  achieved because the proof of \cite[Proposition 2]{Miller} utilizes
  little more than \Souslin's theorem.
\end{propositionproof}

Proposition \ref{maximal:AD} gives rise to an analogous version of
Proposition \ref{coloring}. As every subset of a topological space is
$\aleph_0$-universally \Baire under $\AD$ (see, for example, \cite
[Theorem 7D.2]{Moschovakis}), this can be combined with Proposition
\ref{homomorphism} and \Kanovei's generalization of the $\Gzero$
dichotomy (see \cite{Kanovei}, although the elementary proof of
\cite[Theorem 8]{Miller} can be adapted to obtain the special cases we
need by keeping track of codes as above) to establish Theorem \ref
{main:AD}.

By eliminating the outer induction and the use of \cite[Proposition 2]
{Miller} in the proof of Proposition \ref{maximal}, one obtains a proof
of the weaker result without definability conditions on the sets
involved. Moreover, this proof trivially generalizes to colorings $c \from
X \to \kappa$, for any aleph $\kappa$, and gives rise to an analogous
version of Proposition \ref{coloring}. As a result of \Woodin's ensures
that every subset of an analytic \Hausdorff space is $\kappa$-\Souslin,
for some aleph $\kappa$, under $\ADR$ (see, for example, \cite
[Theorem 32.23]{Kanamori}), this can be combined with Proposition
\ref{homomorphism} and the weakening of \Kanovei's generalization of
the $\Gzero$ dichotomy in which there are no definability constraints on
the coloring (which follows from the simplification of the proof of
\cite[Theorem 8]{Miller} in which the use of \Souslin's theorem is
eliminated) to establish Theorem \ref{main:ADR}.

\bibliographystyle{amsalpha}
\bibliography{bibliography}

\end{document}